\documentclass[psamsfont,a4paper,12pt]{amsart}

\usepackage[english]{babel}   
\numberwithin{equation}{section}

\usepackage[nospace,noadjust]{cite}

\usepackage{setspace}     
\usepackage{enumerate} 	
\usepackage{afterpage} 	
\usepackage[usenames,dvipsnames]{color}
\usepackage{graphicx}	
\renewcommand{\arraystretch}{2}
\usepackage{amsmath} 	
\usepackage{amssymb}	
\usepackage{mathrsfs} 	
\usepackage{amsfonts} 	
\usepackage{latexsym} 	

\usepackage{amsthm} 	
\usepackage{stackrel}       
\usepackage{amscd} 	

\usepackage{tabularx}	
\usepackage{longtable}	

\usepackage{verbatim} 
\usepackage{blindtext} 

\usepackage{amssymb,latexsym,amsmath}    
\usepackage{pstricks,multido,pst-plot}		 
\usepackage{multicol,fancyheadings}			   
\usepackage{float}                       
\usepackage[vcentermath,enableskew]{youngtab}

\usepackage{subfig}

\allowdisplaybreaks

\newcommand{\C}{\mathbb{C}}

\newcommand{\R}{\mathbb{R}}

\newcommand{\Z}{\mathbb{Z}} 
\newcommand{\Q}{\mathbb{Q}}

\def\Tr{\mathop{\rm Tr}}

\renewcommand{\to}{\longrightarrow}

\def\mod{\mathop{\rm{mod}}}

\def\exp{{\mathrm{exp}}}

\newtheorem{theorem}{Theorem} 

\newtheorem*{definition}{Definition}

\newtheorem*{remark}{Remark}

\theoremstyle{definition}

\theoremstyle{remark}

\newtheorem{ind}[]{{\rm\it Indice}}

\usepackage{multicol}

\usepackage{tikz}
\usetikzlibrary{arrows}

\usepackage{paralist}
\usepackage{cite}

\title{From partitions to Hodge numbers of Hilbert schemes of surfaces}

\author{Nate Gillman}
\address{Department of Mathematics \& Computer Science, Wesleyan University, Middletown, CT 06457}
\email{ngillman@wesleyan.edu}

\author{Xavier Gonzalez}
\address{Mathematical Institute, University of Oxford, Oxford, England}
\email{xavier.gonzalez@balliol.ox.ac.uk}

\author[Ono]{Ken Ono}
\address{Department of Mathematics, University of Virginia, Charlottesville, VA 22904}
\email{ken.ono691@gmail.com}

\author[Rolen]{Larry Rolen}
\address{Department of Mathematics,
Vanderbilt University,
Nashville, TN 37240}
\email{larry.rolen@vanderbilt.edu}

\author{Matthew Schoenbauer}
\address{Department of Mathematics, University of Notre Dame, Notre Dame, IN 46556}
\email{mschoenb@nd.edu}

\begin{document}


\begin{abstract}  We celebrate the 100th anniversary of Srinivasa Ramanujan's election as a Fellow
of the Royal Society,  which was largely based on his work with G. H. Hardy on the asymptotic properties of the partition function.
After recalling this  revolutionary work, marking the birth
of the ``circle method'', we present a contemporary example of its legacy in topology.  
We  deduce the equidistribution of Hodge numbers for Hilbert schemes of suitable smooth projective surfaces.

\end{abstract}
\dedicatory{To commemorate the 100th anniversary of Ramanujan's election as a Fellow of the Royal Society}
\maketitle

\section{Introduction}

A {\it partition} is any nonincreasing sequence of positive integers, and the
partition function $p(n)$ counts the number with  size $n$.  Euler established the beautiful
fact that its generating function is given by the infinite product
\begin{equation}\label{GenFcn}
P(q)=\sum_{n=0}^{\infty} p(n)q^n=\prod_{n=1}^{\infty}\frac{1}{1-q^n}=1+q+2q^2+3q^3+5q^4+\dots.
\end{equation}
Ramanujan elegantly made use of this infinite product to prove some of the first deep theorems
about the partition function. Indeed, he used it to prove \cite{Ramanujan}
his well-known congruences
\begin{displaymath}
\begin{split}
p(5n+4)&\equiv 0\pmod{5},\\
p(7n+5)&\equiv 0\pmod{7},\\
p(11n+6)&\equiv 0\pmod{11}.
\end{split}
\end{displaymath}
These congruences have inspired the entire field of partition congruences 
\cite{AhlgrenOno}.

Although partitions are simple to define and the $p(n)$ congruences above are quite beautiful, they turn out to be notoriously difficult to count. 
The following table underscores the nature of this problem by exhibiting the astronomical rate of growth of $p(n)$.

\renewcommand{\arraystretch}{1}
\begin{table}[h] 
\begin{center}
\scalebox{0.9}{\begin{tabular}{|c|c|c|c|c|c|}\hline 
     $n$ & $p(n)$   \\\hline 
     $10$ & $42$  \\\hline 
      $20$ & $627$  \\\hline 
      $40$ & $37,338$  \\\hline 
       $80$ & $15,796,476$  \\\hline 
       $\vdots$ & $\vdots$  \\\hline 

\end{tabular}}
\end{center}
\caption{Values of $p(n)$} \label{table_2}
\end{table}

Hardy and Ramanujan \cite{HardyRamanujan} stunned the mathematical community with their proof of their asymptotic formula.

\begin{theorem}{\text {\rm  (Hardy-Ramanujan, 1918)}}\label{Asymptotic}
As $n\rightarrow +\infty,$ we have
$$
p(n)\sim \frac{1}{4n\sqrt{3}}\cdot e^{\pi \sqrt{2n/3}}.
$$
\end{theorem}

In fact, Hardy and Ramanujan proved a much stronger result than this asymptotic. They showed that $p(n)$ can be approximated by means of a divergent series, sharpening the asymptotic 
in Theorem~\ref{Asymptotic}. Specifically, they showed that there is a sum of similar terms such that for some constant $C$, 
\begin{equation}\label{pofnHRExp}
p(n)=\sum_{j=1}^{C\cdot\sqrt n}E_j(n)+O(n^{-\frac14})
.
\end{equation}
Here, $E_1(n)\sim\frac{1}{4n\sqrt{3}}\cdot e^{\pi \sqrt{2n/3}}$ is the main term, and the later terms have similar asymptotics but for exponentials with smaller multiples of $\sqrt n$. 
However, the series of $E_j(n)$ summed up over all $n$ actually diverges, and so this result falls short of a proper series expansion for $p(n)$.

Twenty years later Rademacher perfected \cite{RademacherPn} this idea and obtained a series expansion which does converge, thus giving an {\it exact} formula for $p(n)$. 
The flavor of both expansions is that they are expressible as sums of Bessel functions times Kloosterman sums. However, Rademacher utilized a different approach which gave different Bessel functions. Although the result is asymptotically the same at each stage, the savings are sufficient to make the sum over all $j$ converge. As we shall see, Rademacher used the same method as Hardy and Ramanujan, namely, their ``circle method,'' though he modified the details in the exact path of integration which led to this improvement. This led to the following formula, where $I_{\frac32}(\cdot)$ is the usual Bessel function, and $A_k(n)$ is the Kloosterman sum
\begin{equation}\label{KloostermanPn}
A_k(n):=\frac12\sqrt{\frac k{12}}\sum_{\substack{d\pmod{24k} \\ d^2\equiv-24n+1\pmod{24k}}}\left(\frac{12}{d}\right)e^{2\pi i \cdot \frac{d}{12k}}
.
\end{equation}

\begin{theorem} {\text {\rm (Rademacher \cite{RademacherPn})}}\label{Exact}
For any natural number $n$, we have
\[
p(n)=\frac{2\pi}{(24n-1)^{\frac34}}\sum_{k\geq1}\frac{A_k(n)}{k}I_{\frac32}\left(\frac{\pi\sqrt{24n-1}}{6k}\right)
.
\]
\end{theorem}

The shape of Rademacher's formula for $p(n)$ would later be understood to arise naturally from the method of {\it Poincar\'e series} by the work of Petersson, Rademacher, and others
(for example, see the exposition in \cite{OurBook}). These are natural modular forms which are built as averages over the translates of suitable special functions under the action of the modular group. In the case of $p(n)$, the generating function is essentially a weight $-1/2$ modular form.

 In the case of half-integral weight Poincar\'e series, formulas such as Rademacher's, understood via the modern theory of Poincar\'e series, can be used to give 
finite, exact formulas for coefficients of modular forms. One of the first important examples of this phenomenon was observed by Zagier
\cite{Zagier} in his work on traces of singular moduli (see also \cite{BringmannOno}).
This idea also applies to $p(n)$. 
Namely,
Rademacher's exact formula for $p(n)$ can be reformulated as a finite sum of values of a single (non-holomorphic) modular function. This fact was first observed by Bringmann and one of the authors
\cite{BringmannOnoPn}, and the phenomenon relies on the fact that (\ref{KloostermanPn}) can be reformulated as a sum over equivalence classes of discriminant $-24n+1$ positive definite integral binary quadratic forms. This observation was refined by Bruinier and one of the authors \cite{BruinierOno} to prove
a much stronger statement.

To make this precise, we let $\eta(\tau):=q^{1/24}\prod_{n=1}^{\infty}(1-q^n)$ ($q:=e^{2\pi i \tau}$ throughout) be Dedekind's weight 1/2 modular form.
Furthermore, we let $E_2(\tau):=1-24\sum_{n=1}^{\infty} \sum_{d\mid n}dq^n$ be the usual weight 2 quasimodular Eisenstein series, and we let
$F(\tau)$ be the weight $-2$ meromorphic modular form
$$
F(\tau):=\frac{1}{2}\cdot \frac{E_2(\tau)-2E_2(2\tau)-3E_2(3\tau)+6E_2(6\tau)}{\eta(\tau)^2
\eta(2\tau)^2 \eta(3\tau)^2\eta(6\tau)^2}
=q^{-1}-10-29q-\dots.
$$
Using the convention that $\tau=x+iy$, with $x, y\in \R$,  we define the weight 2 weak Maass form
\begin{equation}
\mathcal{P}(\tau):=-\left(\frac{1}{2\pi i}\cdot \frac{d}{d\tau}+\frac{1}{2\pi y}\right) F(\tau)=
\left(1-\frac{1}{2\pi y}\right)q^{-1}+\frac{5}{\pi y}+\left(29+\frac{29}{2\pi y}\right)q+\dots.
\end{equation}

The finite algebraic formula for $p(n)$ is given in terms of the {\it singular moduli} for $\mathcal{P}(\tau)$, the values of this weak Maass forms
at CM points. More precisely, we use discriminant $-24n+1=b^2-4ac$ positive definite integral binary quadratic forms
$$
Q(x,y)=ax^2+bxy+cy^2,
$$
with the property that $6\mid a$. The congruence  subgroup $\Gamma_0(6)$ acts on these forms, and we let $\mathcal{Q}_n$ be the (finitely many) equivalence
classes with $a>0$ and $b\equiv 1\pmod{12}$. If $Q(x,y)$ is such a form, then we let $\alpha_Q$ be the unique point in the upper-half
of the complex plane for which $Q(\alpha_Q,1)=0$. By the theory of complex multiplication, these values are algebraic, and they generate
ring class field extensions of $\Q(\sqrt{-24n+1})$. We then define their trace by
\begin{equation}
\Tr(n):=\sum_{Q\in \mathcal{Q}_n} \mathcal{P}(\alpha_Q).
\end{equation}
In terms of this notation, we have the following pleasing theorem.

\begin{theorem}{\text {\rm (Bruinier-Ono \cite{BruinierOno}, 2013)}}\label{Finite}
If $n$ is a positive integer, then we have
$$
p(n)=\frac{1}{24n-1}\cdot \Tr(n).
$$
The numbers $\mathcal{P}(\alpha_Q)$, as $Q(x,y)$ varies over the finitely many classes in $\mathcal{Q}_n$, form a multiset of algebraic numbers which is the union
of Galois orbits for the discriminant $-24n+1$ ring class field. Moreover,  for each $Q\in \mathcal{Q}_n$ we have that
$(24n-1)\mathcal{P}(\alpha_Q)$ is an algebraic integer.
\end{theorem}

\begin{remark}
Larson and one of the authors {\text {\rm \cite{LarsonRolen}}} established the precise integrality properties of the values
$\mathcal{P}(\alpha_Q)$.
\end{remark}

The proofs of Theorems~\ref{Asymptotic}, \ref{Exact} and \ref{Finite} depend critically on the fact that
the generating function in (\ref{GenFcn}), where $q:=e^{2\pi i \tau}$, satisfies
$q^{-1/24}P(q)=1/\eta(\tau)$. 
It is now well understood that the circle
method can be applied to modular forms whose poles are supported at cusps. Moreover, for those modular forms
that have non-positive weight, the method typically offers exact formulas as in Theorem~\ref{Exact}.
Here we describe recent developments
in topology
in which such modular forms that can be written as infinite products that resemble \eqref{GenFcn} arise as  generating functions of topological invariants. Thus, by applying the circle method of Hardy and Ramanujan as perfected by Rademacher, we obtain exact formulas that provide insight into the distribution of these topolological invariants.

To discuss this connection between topology and Ramanujan's legacy, we recall the foundational importance of \emph{topological invariants}. One of the broad goals of topology  is to determine whether two particular spaces have the same topological, differentiable, or complex analytic structure. 
When this is the case, one can often find an isomorphism between the two spaces, identifying them in a way that respects this structure. 
It is, however, generally more difficult to prove that two spaces are fundamentally distinct. 
\textit{Topological invariants} assign numbers, groups, or other mathematical objects to spaces in such a way that isomorphic spaces yield the same output.
In this way, invariants are useful for distinguishing dissimilar spaces.

Here the spaces we are concerned with are complex manifolds. 
An important class of invariants known as the \textit{Hodge numbers} $h^{s,t}$ belong to manifolds of this type. For any $n$-dimensional complex manifold $M$ and any $0 \leq s,t,\leq n$, the Hodge number $h^{s,t}(M)$ gives the dimension of a certain vector space of differential forms on $M$. 
For the manifolds we will be concerned with, important topological invariants such as the Betti numbers and signature arise as linear combinations of the Hodge numbers
(see \cite{Wells}).

There are many ways of constructing new spaces from old, and when we study topology we want to understand how invariants interact with these constructions. In algebraic geometry, the $n^{th}$  \textit{Hilbert scheme} of a projective variety $S$ is a projective variety $\mathrm{Hilb}^n(S)$ that can be thought of as a smoothed version of the $n^{th}$ symmetric product of $S$ (for example, see \cite{Hilbert_schemes}).
The $n$-th symmetric product of a manifold $M$ admits a simple combinatorial interpretation: outside of a negligible subset, the symmetric product is the collection of subsets of $M$ of size $n$ assembled as a manifold in its own right.
Interestingly, the Hodge numbers of a complex projective surface $S$ determine the Hodge numbers of $\mathrm{Hilb}^n(S)$ for all $S$
in a very pleasing combinatorial way.
 This statement is captured in the following beautiful theorem of G\"ottsche \cite{Gottsche}

\begin{theorem}[G\"ottsche]\label{thm:gottsche} If $S$ is a smooth projective complex surface, then we have that
\begin{equation} \sum_{\substack{n\geq 0\\\mathbf{0\leq s,t \leq 2n}}}(-1)^{s+t}h^{s,t}(\mathrm{Hilb}^n(S))x^{s-n}y^{t-n}q^n =\prod_{n =1}^\infty \frac{\prod_{s +t \mathrm{\ odd}} (1- x^{s-1}y^{t-1}q^n)^{h^{s,t}(S)}}{\prod_{s +t \mathrm{\ even}} (1- x^{s-1}y^{t-1}q^n)^{h^{s,t}(S)}}. \label{eq:goettsche} \end{equation}
\end{theorem}

The fortunate feature of this formula is that the Hodge numbers $h^{s,t}(\mathrm{Hilb}^n(S))$ are  prescribed by the  infinite product in (\ref{eq:goettsche}) which can be specialized to obtain modular forms. 
In these cases we will use the circle method to obtain exact formulas, as well as asymptotic and distributional information, for these Hodge numbers for a certain class of complex projective surfaces. This work generalizes previous work by some of the authors \cite{reu_paper}.

We pursue this task in the same spirit that led Ramanujan to bring forward new information about the partition numbers using the modularity of $\eta(\tau)$. 
Indeed, the process of taking symmetric powers of surfaces is inherently combinatorial, and in the spirit of Ramanujan gives rise precisely to an infinite family of topologically inspired ``partition problems'' encoded in (\ref{eq:goettsche}).
Hence, we expect the circle method to apply.
In other words, the role of infinite products in partition theory offers a glimpse of G\"ottsche's theorem as a device which mirrors the assembly required to build $\mathrm{Hilb}^n(S)$.
Although the combinatorial object we are studying, $\mathrm{Hilb}^n(S)$, arises in a different field of mathematics than partitions, we find that it is Ramanujan's insight and his novel use of modular forms that illuminates the path to obtaining new and interesting information about these spaces.

Towards this end, we begin by collecting the Hodge numbers for a complex manifold $M$ in a generating function called the \emph{Hodge polynomial}:
\begin{equation*}
     \chi_{\mathrm{Hodge}}(M)(x,y) := x^{-d/2}y^{-d/2}\sum_{s,t}h^{s,t}(M)(-x)^s(-y)^t.
\end{equation*} Henceforth we refer to the generating function for the Hodge polynomials of Hilbert schemes of a smooth projective complex surface $S$ as
\begin{equation}\label{eq:def_z}
    Z_S(x,y;\tau):=\sum_{n = 0}^\infty \chi_{\mathrm{Hodge}}(\mathrm{Hilb}^n(S))(x,y)q^n. 
\end{equation}By specializing \eqref{eq:def_z} appropriately, we obtain generating functions for a variety of topological invariants (see \cite{reu_paper}).
In order to study the distributional properties of Hodge numbers, we consider
\begin{equation}\label{eq:def_gamma}
    \gamma_{S}(r_1,\ell_1, r_2, \ell_2;n):= \sum_{\substack{{t \equiv r_1} \mod{\ell_1}\\ {s \equiv r_2 } \mod{\ell_2} }} (-1)^{s+t}h^{s+n,t+n}(\mathrm{Hilb}^n(S)),
\end{equation}
and compile the generating function
\begin{equation}\label{eq:def_C}
     C_S(r_1,\ell_1, r_2, \ell_2;\tau) := \sum_{n \geq 0} \gamma_S(r_1,\ell_1, r_2, \ell_2;n) q^n.
\end{equation}
We would like to determine when the Hodge numbers of a surface $S$ are equidistributed.
We define such an \emph{equidistribution} as follows:

\begin{definition}
Let $S$ be a smooth projective complex surface. We say that $S$ has $(\ell_1,\ell_2)$-equidistribution if for some $\mathcal{R} \subseteq \Z/{\ell_1}\Z \times \Z/{\ell_2}\Z$ we have, as $n\to\infty$, 
$$ \gamma_S(r_1,\ell_1,r_2,\ell_2;n) \sim \gamma_S(r_1',\ell_1,r_2',\ell_2;n)
$$
for all $(r_1,r_2),(r_1',r_2') \in \mathcal{R}$, and 
$$\gamma_S(r_1,\ell_1,r_2,\ell_2;n) = 0$$
for all $(r_1,r_2) \not\in \mathcal{R}$ and all $n>0$. 
\end{definition}
Since the generating function in \eqref{eq:def_C} arises as a linear combination of specializations of 
\eqref{eq:def_z} according to
\begin{equation}
\label{eq:roots}
 C_S(r_1,\ell_1, r_2, \ell_2;\tau) = \frac{1}{\ell_1\ell_2} \sum_{\substack{ j_1 \mod{\ell_1} \\ j_2\mod{\ell_2}}}  \zeta^{-j_2r_2}_{\ell_2} \zeta^{-j_1r_1}_{\ell_1} Z_S(\zeta_{\ell_2}^{j_2}, \zeta^{j_1}_{\ell_1};\tau),
\end{equation}
where $\zeta_\ell$ is a primitive $\ell^{th}$ root of unity, determining the behavior of specializations of \eqref{eq:def_z} is useful for studying the distribution of Hodge numbers. We can express these functions $$
  Z_S(\zeta_{\ell_1}^{r_2}, \zeta^{r_2}_{\ell_2};\tau)  = : \sum _{n \geq 0}  \xi_S(r_1,\ell_1,r_2,\ell_2;n) q^n
$$ in terms of $\eta(\tau)$ and \emph{generalized Dedekind} $\eta$ \emph{functions}, which are defined on page 187 of \cite{Schoeneberg} as $$ \eta_{\left(u,v,N\right)}\left(\tau\right): = \alpha_N\left(u,v\right)e^{\pi i P_2\left(u/N\right)\tau} \prod_{\substack{m>0 \\ m \equiv u \mod{N}}}\left(1- \zeta_N^{v}e^{2 \pi i \tau m/N}\right) \prod_{\substack{m>0 \\ m \equiv -u \mod{N}}}\left(1- \zeta_N^{-v}e^{2 \pi i \tau m/N}\right)\mathbf{,}$$ {where $\alpha_N(u,v)$ is given by }
$$
\alpha_N(u,v):=
\begin{cases}
\left(1 - \zeta_N^{-v}\right)e^{\pi i P_1\left( \frac{v}{N}\right)} & u \equiv 0 \text{ and } v \not \equiv 0\  \pmod N, \\
\hfil 1  & \text{{otherwise.}}
\end{cases}
$$
  Here, \textbf{$P_1(x)=\{x\}-1/2$} and $P_2(x):= \{x\}^2 -\{x\} + 1/6$. By (\cite{Schoeneberg}, pg. 200) we have that for $u,v,N \in \Z$, $(u, v ) \not \equiv 0 \mod N$, $\eta_{u,v,N}^{N_1}(\tau)$ is a modular function on $\Gamma(N)$, where $N_1 = 12 N^2 / \gcd(6,N)$. The explicit transformation law for $\eta_{(u,v,N)}(\tau)$ on all of $\mathrm{SL}_2(\mathbb{Z})$ is shown on page 198 of \cite{Schoeneberg}. We exploit this transformation law in Section \ref{ssc:prf_sketch} to find an exact formula for the $\xi_S(r_1,\ell_1,r_2,\ell_2;n)$ for suitable surfaces.
To make this precise, suppose that $M$ is a $d$-dimensional complex manifold. We let $\chi(M)$ denote its
Euler characteristic, and we let
$\sigma(M)$ denote the signature of the intersection pairing on $H^d(M)$. Then we obtain the
following exact formulas, 
{where
$L:=\mathrm{lcm(\ell_1,\ell_2)}$, 
$H:=H(\iota_2)$ is given by \eqref{eq:Hdefn},
$a_j$ is a Fourier coefficient defined in \eqref{eq:Zstar},}
$B_k$ is a Kloosterman sum defined in \eqref{eq:Kloosterman_def}, 
and $I^*$ is a scaled modified Bessel function of the first kind defined in \eqref{eq:I_def}.

\begin{theorem}\label{thm:exact_formulas}
Let $S$ be a smooth projective complex surface such that $\chi(S) \geq 0$ and $\chi(S) \geq \sigma(S)$.
Then,
 \begin{equation} \xi_S(r_1,\ell_1,r_2,\ell_2;n)=
 2\pi \alpha \sum_{\substack{\iota_1 \mod{L} \\ \iota_2 \mod{L}}} 
    \sum_{j<-L{H}} 
    \sum_{\substack{k=1 \\ k \equiv \iota_2 \mod{L}}}^\infty 
    \frac{\alpha'a_j}{k^G}B_k(j,L,\iota_1;n)
    I^*(\iota_1,\iota_2,j,k;n), \label{eq:thm5}. \end{equation}
 \end{theorem}
 
  \begin{remark}
  Theorem \ref{thm:exact_formulas} holds for a large class of surfaces $S$. In particular, it holds for all but finitely many Hodge structures in each birational equivalence class. In addition, the Enriques-Kodaira  Classification  Theorem implies that for surfaces of non-general type, the only minimal models which do not satisfy the hypotheses of Theorem \ref{thm:exact_formulas} are ruled surfaces of genus $g$, where $g \geq 2$ (see \cite{reu_paper}).
 \end{remark}
 \begin{remark}
  By the following equality,
 $$\gamma_{S}(t-n,2n+ 1, s-n, 2n+1;n)=  (-1)^{s+t}h^{s,t}(\mathrm{Hilb}^n(S)),$$
 Theorem \ref{thm:exact_formulas} and Equation (\ref{eq:roots}) together give exact formulas for  $h^{s,t}(\mathrm{Hilb}^n(S))$ for surfaces $S$ such that $\chi(S) \geq 0$ and $\chi(S) \geq \sigma(S)$. 
 \end{remark}
 
Building on work in \cite{manschot_rolon}, from which it follows that for $S$ a K3 surface $\gamma_S(r, \ell, 0, 1; n) \sim \gamma_S(r', \ell, 0, 1; n)$ as $n \to \infty$, and \cite{reu_paper}, which described equidistribution in the case $\ell_1 = \ell_2 = 2$, we use the asymptotics derived from the exact formula in Theorem \ref{thm:exact_formulas} to make the following statement about the equidistribution of Hodge numbers for appropriate surfaces:
 
\begin{theorem}\label{thm:equidistribution}
Let $S$ be a smooth projective complex surface such that $\chi(S) \geq \sigma(S)$. Then $S$ has $(\ell_1,\ell_2)$-equidistribution if and only if one or more of the following holds:
\begin{enumerate}
\item 
$h^{1,0} = 0$, $h^{2,0}=0$, $\mathcal{R} = \{ (r_1,r_2) \ | \ r_1  \equiv r_2 \mod \gcd(\ell_1,\ell_2)\}$,
\item
$h^{1,0} = 0$, $h^{2,0} >0$, $\mathcal{R} = \{ (r_1,r_2) \ | \ r_1  \equiv r_2 \mod \gcd(\ell_1,\ell_2,2)\}$,
\item $\chi(S) + \sigma(S) = 0$, $\chi(S) \neq 0$, { $\min\{\ell_1,\ell_2\} = 1$,} $\mathcal{R} = \{(0,0)\}$,
\item $\chi(S) + \sigma(S) = 0$, $\chi(S) = 0$, { $\min\{\ell_1,\ell_2\} = 1$,} $\mathcal{R} = \emptyset$,
{ \item  $\chi(S) \neq 0$, $\ell_1,\ell_2=1$, $\mathcal{R} = \{(0,0)\}$,
\item $\chi(S) = 0$, $\ell_1,\ell_2=1$, $\mathcal{R} = \emptyset$,}
\item $h^{1,0} >0$, $\chi(S) + \sigma(S) >0$, $\mathcal{R} =  \Z/{\ell_1}\Z \times \Z/{\ell_2}\Z$, and $ { \Lambda(}0,0) < { \Lambda(}j_1/\ell_1,j_2/\ell_2)$, for all $(j_1,j_2) \not \equiv (0,0)$, where $$ { \Lambda(}x,y) :=  h^{1,0} \left( P_2\left(x\right) + P_2\left(y\right)\right) 
    - h^{0,0} P_2\left( x+y\right) 
    - h^{2,0} P_2\left( x-y\right)
$$. 
\end{enumerate}
Case (7) occurs whenever $\min\{\ell_1,\ell_2\} = 1$, and holds for only finitely many $(\ell_1,\ell_2)$ such that $\min\{\ell_1,\ell_2\}> 1$. In addition, case (7) only occurs when $\gcd(\ell_1,\ell_2)=1$.
\end{theorem}

\section{The partition function}

Here we describe the history and main idea of the circle method for the partition function. {This will use the modularity of the generating function for the partition function. We should point out that though the circle method is especially convenient to apply in such cases, many important applications of the circle method are possible in non-modular situations, for example as explored by Hardy and Littlewood
\cite{Vaughan}  in their study of Waring's problem.} The basic idea is simple. Recall from \eqref{GenFcn} that the generating function for $p(n)$ satisfies a product form first discovered by Euler:
\[
P(q)=\prod_{n=1}^{\infty}\frac{1}{1-q^n}
.
\]
Cauchy's residue theorem can then be used to isolate any of the coefficients of this expansion. Specifically, if we divide $P(q)$ by $q^{n+1}$, then as a function of $q$, $P(q)/q^{n+1}$ has residue $p(n)$, and so 
\[
p(n)=\frac1{2\pi i}\int_C\frac{P(q)}{q^{n+1}}dq
,
\]
where $C$ is any simply closed path around the origin contained in the unit circle, traversed counterclockwise. Though this idea is straightforward in principle,
great care must be taken in choosing a suitable path $C$ so that the integral can be closely estimated. To determine a ``good'' path,
one must first consider the location of the poles of $P(q)$. Again, this is furnished by Euler's product formula \eqref{GenFcn}, which shows that the poles lie exactly at the roots of unity. This justifies the earlier claim that we may take any path which doesn't cross this wall of singularities, and gives a first indication of how to estimate this integral. 
We would like to split this up into a ``main term'' and an error term, and so it is important to study where most of the contribution of the integral will come from. If we choose a path which 
approaches roots of unity quite closely, then the majority of the integral will come from the parts of the path near these poles. However, not all poles will contribute equally to the size of the integral. 
An analysis of the generating function $P(q)$ shows that ``near'' primitive $j$-th order roots of unity, the size is much smaller the larger $j$ is. Thus, the main contribution is from the pole at $q=1$, secondary terms come from the pole at $q=-1$, and the third order of contribution comes from the third order roots of unity. In fact, the explanation of the terms in the expansion \eqref{pofnHRExp} is that the function $E_j(n)$ is an 
approximation of the behavior of a suitable Cauchy integral near the primitive $j$-th order roots of unity. As an illustration of this numerical phenomenon, consider the following table of values of $P(q)$ near roots of unity, where the columns correspond to the different roots of unity $\zeta$, and where the values of $t$ in the rows correspond to evaluating $|P(\zeta\cdot e^{-t})|$.

\begin{table}[h] 
\small
\begin{center}
\scalebox{0.9}{\begin{tabular}{|c|c|c|c|c|}\hline 
  &   $\zeta=1$&$\zeta=-1$  &$\zeta=-\frac12\pm\frac{\sqrt{-3}}{2}$  &$\pm i$  \\\hline 
   t=0.5 & 7.4& 0.87&0.68 &0.66 \\ \hline
   t=0.3 & 51.3& 1.2&0.68 &0.60 \\ \hline
   t=0.1 & $1.7\cdot 10^6$&10.8 &1.3 &0.70  \\ \hline
   t=0.01 & $1.1\cdot 10^{70}$&$4.1\cdot 10^{16}$ &$6.0\cdot 10^6$ & 2325.4
    \\ \hline
\end{tabular}}
\end{center}
\caption{Illustrative values of $P(q)$ near roots of unity} \label{table_eta_values}
\end{table}
The exact description of $P(q)$ near roots of unity is afforded by
the modular transformation properties of the Dedekind-eta function,
for which
$$
P(q)=\frac{q^{\frac{1}{24}}}{\eta(\tau)}
.
$$
Essentially, $\eta(\tau)$ is a level one modular form of weight $1/2$, with a multiplier system consisting of $24$-th roots of unity. 
That is, for any $\gamma=\left(\begin{smallmatrix}a & b\\ c& d\end{smallmatrix}\right)\in\operatorname{SL}_2(\mathbb Z),$
\[
\eta(\gamma\cdot\tau)=\omega_\gamma(c\tau+d)^{\frac12}\eta(\tau),
\]
where $\omega_\gamma^{24}=1$.
Specifically, the numbers $\omega_\gamma$ are determined by the values at the matrices
\[
T:=\begin{pmatrix} 1 & 1\\ 0 &  1\end{pmatrix}, \quad\quad\quad S:=\begin{pmatrix}0 & -1\\ 1 & 0\end{pmatrix}
\]
as follows:
\[
\eta(\tau+1)=e^{\frac{2\pi i}{24}}\eta(\tau),
\quad\quad\quad
\eta\left(-\frac{1}{\tau}\right)=\sqrt{-i\tau}\eta(\tau)
.
\]
Using matrices to ``connect'' any root of unity to the point at infinity, and dealing with the elementary factor $q^{\frac{1}{24}}$ yields the desired expansions near the cusps. 
For instance, as $q$ approaches $1$ radially from within the unit disk, $q=e^{2\pi i \tau}$ with $\tau$ tending to zero, we have
\[P(q)\sim\sqrt{-i\tau}e^{\frac{\pi i}{12\tau}}.\]
Such calculations suffice to estimate the values of $P(q)$ along any desired path. Now, we must discuss which path one should choose. We will follow Rademacher's choice, which yields his exact formula above. The first natural choice, as suggested by the name of the method, could be to let $C$ be a circle centered at the origin. In the upper half plane, that is as a function of $\tau$, this corresponds to choosing a horizontal path with endpoints 1 unit apart. Rademacher replaced this path by an increasingly large number of mutually tangent circles, which are the well-known {\it Ford circles}. For each rational number, there is a single Ford circle tangent to it. For each cutoff $N$, Rademacher considered the subset of these circles which are tangent to the rational numbers with denominator less than or equal to $N$. The exact description is beautifully expressed using the theory of Farey fractions. Explicitly, these circles are precisely the image of the line $\operatorname{Im}(\tau)=1$ in the upper half plane under the action of $\operatorname{SL}_2(\mathbb Z)$.

The image on the following page depicts the Ford circles centered around fractions in $[0,1]$ with denominator at most $N=4$.
For each $N$, Rademacher's path is to traverse the arcs of these circles, starting at $i$, moving along each arc until the next included Ford circle is intersected, until finally arriving at $i+1$. In the image for $N=4$, the path travels along the solid circular arcs. Rademacher then expressed the Cauchy integral as a sum over each of these paths, and, under a suitable change of variables, then used the modular properties of the eta function to expand the values of $P(q)$ on each of these arcs. These expansions of $P(q)$ can then naturally be split up into a main term and an error term, the main features being that the main terms become integrals which can be evaluated exactly as Bessel functions, and the error terms can be explicitly bounded.

After the work of Rademacher, it became apparent that Rademacher's series were in fact Poincar\'e series in disguise. In particular, there are special functions which, when averaged under the slash operator over the modular group, have Fourier expansions which give the same expansions. As such Poincar\'e series are always a basis of the space of all weakly holomorphic modular forms up to cusp forms (as one can explicitly match them to any principal part and constant term which always determines a modular form up to a cusp form). For forms of negative weight such as $1/\eta(\tau)$ this procedure always yields an exact formula. More details on the method of Poincar\'e series and how it can be applied to general weakly holomorphic and mock modular forms can be found in \cite{OurBook}.

\begin{figure}
\resizebox{13.0cm}{!}{\begin{tikzpicture}
    \begin{scope}[thick,font=\scriptsize]

    \draw [->] (-0.5,0) -- (11,0) node [above left]  {$\Re\{z\}$};
    \draw [->] (0,-0.5) -- (0,11) node [below right] {$\Im\{z\}$};

    \draw [thin] (10,-0.25) -- (10, 0.25);
    \node at (10,-0.5) {$1$};

    \draw [thin] (-0.25, 10) -- (0.25, 10);
    \node at (-0.5, 10) {$i$};

    \draw [thin] (5,-0.15) -- (5, 0.15);
    \node at (5,-0.5) {$1/2$};

    \draw [thin] (10/3,-0.15) -- (10/3, 0.15);
    \node at (10/3,-0.5) {$1/3$};

    \draw [thin] (20/3,-0.15) -- (20/3, 0.15);
    \node at (20/3,-0.5) {$2/3$};

    \draw [thin] (2.5,-0.15) -- (2.5, 0.15);
    \node at (2.5,-0.5) {$1/4$};

    \draw [thin] (7.5,-0.15) -- (7.5, 0.15);
    \node at (7.5,-0.5) {$3/4$};
    

    \draw [blue,domain=-62:90] plot ({0+5*cos(\x)}, {5+5*sin(\x)});
    \draw [dotted,blue,domain=-90:-62] plot ({0+5*cos(\x)}, {5+5*sin(\x)});
    \draw [blue,domain=90:242] plot ({10+5*cos(\x)}, {5+5*sin(\x)});
    \draw [dotted,blue,domain=242:270] plot ({10+5*cos(\x)}, {5+5*sin(\x)});

    \draw [blue,domain=-23:203] plot ({5+1.25*cos(\x)},{1.25+1.25*sin(\x)});
    \draw [dotted,blue,domain=203:337] plot ({5+1.25*cos(\x)}, {1.25+1.25*sin(\x)});

    \draw [blue,domain=23:190] plot ({10/3+(5/9)*cos(\x)},{(5/9)+(5/9)*sin(\x)});
    \draw [dotted,blue,domain=127:383] plot ({10/3+(5/9)*cos(\x)},{(5/9)+(5/9)*sin(\x)});
    \draw [blue,domain=-10:157] plot ({20/3+(5/9)*cos(\x)},{(5/9)+(5/9)*sin(\x)});
    \draw [dotted,blue,domain=157:350] plot ({20/3+(5/9)*cos(\x)},{(5/9)+(5/9)*sin(\x)});

    \draw [blue,domain=10:120] plot ({10/4+(10/32)*cos(\x)},{(10/32)+(10/32)*sin(\x)});
    \draw [dotted,blue,domain=120:370] plot ({10/4+(10/32)*cos(\x)},{(10/32)+(10/32)*sin(\x)});
    \draw [blue,domain=60:170] plot ({30/4+(10/32)*cos(\x)},{(10/32)+(10/32)*sin(\x)});
    \draw [dotted,blue,domain=170:420] plot ({30/4+(10/32)*cos(\x)},{(10/32)+(10/32)*sin(\x)});
    \end{scope}
\end{tikzpicture}}
    \caption{Ford Circles for $N=4$}
    \end{figure}
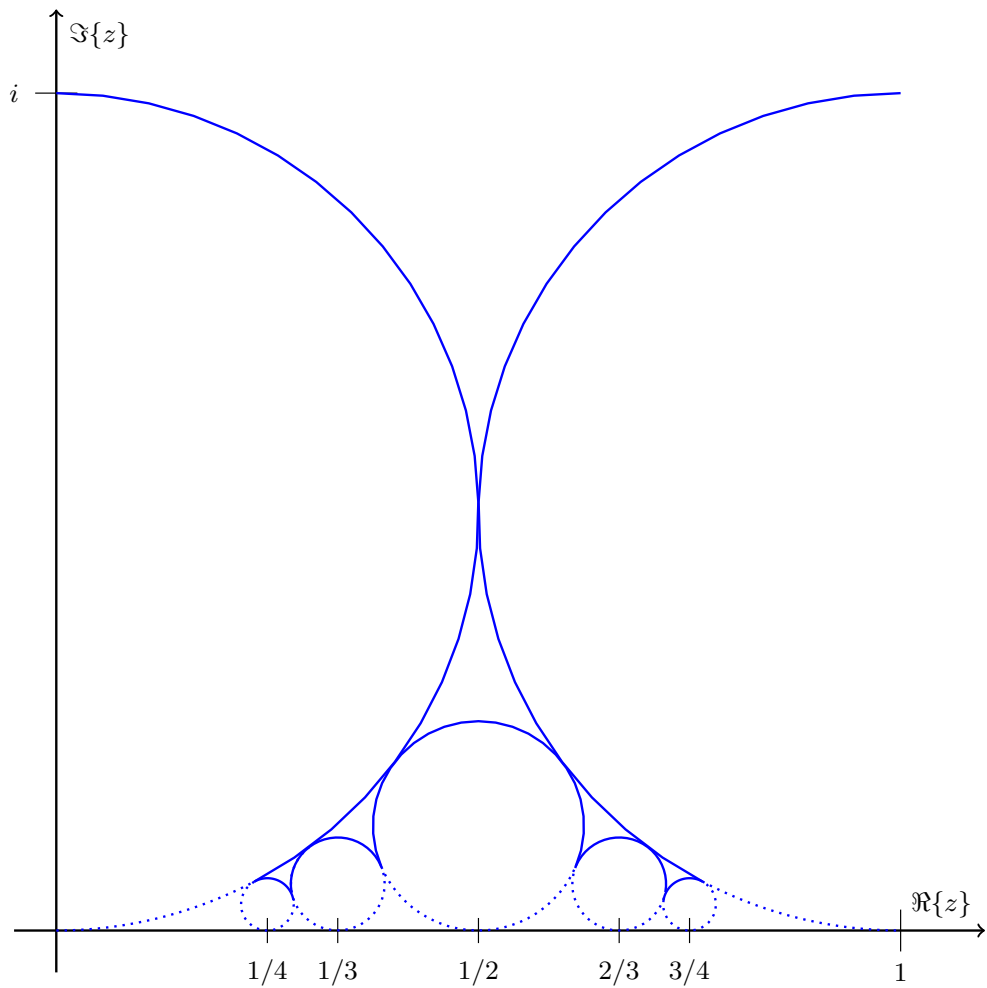

\section{Hodge numbers for Hilbert schemes of surfaces}

We now apply the circle method to obtain exact formulas for Hodge numbers for Hilbert schemes. In the next two subsections we sketch 
the proofs of Theorems~\ref{thm:exact_formulas} and  \ref{thm:equidistribution}, which generalize
earlier work by the authors contained in \cite{reu_paper}.
In the last subsection, we illustrate these results in the case of $S=\C\mathbb{P}^2$.

\subsection{Sketch of the proof of Theorem \ref{thm:exact_formulas}: Exact formulas}\label{ssc:prf_sketch}

For the sake of simplicity, we  will only consider the case $r_2 = 0$, $\ell_2=1$, as the case $r_2 \neq 0$ follows \textit{mutatis mutandis}. 

\begin{proof}[Sketch of the proof of Theorem~\ref{thm:exact_formulas}]
By Cauchy's {  residue theorem}, we obtain the following integral expression for the coefficient of $q^n$ inside $Z_S(\zeta_\ell^r,1;n)$,
\[\xi_S(r,\ell,0,1;n)=\frac{1}{2\pi i}\int_C\frac{Z_S(\zeta_{\ell}^r,1;q)}{q^{n+1}}dq,\]
where we choose $C$ to be a circle of radius $e^{-2\pi N^{-2}}$.

While Rademacher estimated his contour integral by decomposing the line segment in the upper half plane into arcs of Ford circles, we decompose
$C$ into Farey arcs $\Xi_{h,k}$ as in \cite{Rademacher}, yielding
\[\xi_S(r,\ell,0,1;n)=\sum_{\substack{0\leq h<k\leq N\\\text{gcd}(h,k)=1}}\frac{1}{2\pi i}\int_{\Xi_{h,k}}\frac{Z_S\left(\zeta_{\ell}^{r},1;q\right)}{q^{n+1}}dq.\]
We note that with this path of integration, we recover Rademacher's formula for $p(n)$ in Theorem \ref{Exact}. Both methods make use of the same transformation formula and approximate the contribution from the cusps representatives with the same Bessel functions, but differ in the process of bounding errors.
Since we can express our generating function as the modular form
\begin{equation}
    \label{eq:Z_specialized}
 Z_S(\zeta_{\ell}^{r},1;\tau)=
\dfrac{ \alpha q^{\frac{\chi\left(S\right)}{24}} }{\eta_{\left(0,r, \ell\right)}\left(\tau\right)^{(\chi(S) + \sigma(S))/4 } \eta \left(\tau\right)^{(\chi(S) -\sigma(S))/2} }, 
\end{equation} 
where $\alpha$ depends only on $r$ and  $\ell$,
by \cite{Schoeneberg} we have for $(h,k)=1$, $hh' \equiv -1 \  \mod{k}$, and $\mathrm{Re}(z) >0$,
$$
Z_S\left(\zeta_{\ell}^{r}, 1; \frac{iz +h}{k}\right) 
    =
    \omega (h,k) \alpha  \alpha'   \cdot z^{-G} \cdot  \exp\left(-\frac{2 \pi}{k}  \left( \frac{\chi\left(S\right)}{24}z + \frac{H}{z}\right) \right)
     \cdot    Z^*\left(\frac{iz ^{-1} + h'}{k} \right),
$$
where
$$\omega(h,k) := \exp\left( -\pi i/4\cdot  \left(2\left(\chi(S)-\sigma(S)\right) \cdot s(h,k) + 
\left(\chi(S)+\sigma(S)\right) \cdot s_{{ (r,\ell)}}(h,k)
\right)
\right)
$$
is a root of unity built from the following generalized Dedekind sum,
$$s_{{ (r,\ell)}}(h,k) = \sum_{\lambda \mod{k}}\left(\left( \frac{ \lambda }{k}\right)\right)\left(\left( \frac{h\lambda}{k} + \frac{  r}{\ell}\right)\right),
$$
$s(h,k)=s_{(0,1)}(h,k)$, the constant $\alpha'$ depends on $h$ and $k$ mod $\ell$,  $H{:=H(k)}$ is the order of the zero at the cusp $h/k$, {given by}
\begin{equation}\label{eq:Hdefn}
\frac{1}{2} 
    \left(h^{1,0} \left( P_2\left(\frac{k  r_1}{\ell_1}\right) + P_2\left(\frac{k  r_2}{\ell_2}\right)\right) 
    - h^{0,0} P_2\left(\iota_2 \left( \frac{r_1}{\ell_1} +  \frac{r_2}{\ell_2}\right)\right) 
    - h^{2,0} P_2\left(\iota_2 \left( \frac{r_1}{\ell_1} -  \frac{r_2}{\ell_2}\right)\right) - \frac{h^{1,1}}{12}\right),
\end{equation}
$G$ is the weight of $Z_S$, and 
\begin{equation}\label{eq:Zstar}
Z^*(\tau)=\sum_{j\geq 0}a_je^{2\pi i\tau j/\ell}.
\end{equation} 
{The modular transformation law for generalized Dedekind $\eta$ functions found in Chapter 8 of \cite{Schoeneberg} allows us to write $Z^*(\tau)$ as a quotient of infinite products. Using this, one can calculate $a_j$ via a simple product expansion.
}

Since $Z_S(\zeta_\ell^r, 1; \tau)$ is modular on $\Gamma(\ell)$, we deal with the multiple cusps in the spirit of Poincar\'e series by summing over all the Farey arcs near a given cusp, defining
\[S(\iota_1, \iota_2,N;n) := \sum_{\substack{0\leq h<k\leq N\\ \mathrm{gcd}(h,k)=1 \\(h,k) \equiv (\iota_1, \iota_2) { \ \mod{\ell}}}}\frac{1}{2\pi i}\int_{\Xi_{h,k}}\frac{Z_S\left(\zeta_{\ell}^{r},1;q\right)}{q^{n+1}}dq,\]
from which it immediately follows that
\begin{equation}\label{eq:stratify_s}
    { \xi_S(r,\ell,0,1;n)}= \sum_{\substack{\iota_1 \mod{\ell} \\ \iota_2 \mod{\ell}}} S(\iota_1, \iota_2,N;n).
\end{equation}
Once we apply the transformation formula, substitute the series expansion for $Z^*((iz^{-1}+h')/k)$, and then make the variable transformation $w=N^{-2}-i\theta$, we obtain
\begin{align*}
    S(\iota_1, \iota_2,N;n) = 
    & \alpha\alpha'\sum_{j=0}^\infty 
    \sum_{\substack{k=1 \\ k \equiv \iota_2 \mod{\ell}}}^N
    \sum_{\substack{0\leq h<k\\ \mathrm{gcd}(h,k)=1 \\h \equiv \iota_1 \mod{\ell}}}
    B_{h,k}(j,\ell;n)
    \\
    &\cdot \frac{a_j }{k^G}
    \int_{\vartheta_{h,k}'}^{\vartheta_{h,k}''}w^{-G} 
    \exp\left[
    w\left(2\pi n-\frac{\pi\chi(S)}{12}\right)
    +\frac{1}{w}\left(-\frac{2\pi H}{k^2}-\frac{2\pi j}{k^2\ell}\right)
    \right]
     d\theta.
\end{align*} %
Here,
\begin{equation}
  \label{eq:Kloosterman_def}
    B_k(j,\ell,\iota_1;n)
    :=\sum_{\substack{0\leq h<k\\ \mathrm{gcd}(h,k)=1 \\h \equiv \iota_1 \mod{\ell}}}B_{h,k}(j,\ell,n)
    :=\sum_{\substack{0\leq h<k\\ \mathrm{gcd}(h,k)=1 \\h \equiv \iota_1 \mod{\ell}}}
    \omega(h,k)\cdot\exp\left[-\frac{2\pi inh}{k}+\frac{2\pi ih'j}{k\ell}\right]
\end{equation}
is a Kloosterman sum.

Integrating as in \cite{reu_paper}, we find that
\begin{equation}\label{eq:exact_formula_S}
    S(\iota_1, \iota_2,N;n)=2\pi \alpha\alpha' 
    \sum_{j<-\ell H} 
    \sum_{\substack{k=1 \\ k \equiv \iota_2 \mod{\ell}}}^N
    \frac{B_k(j,\ell,\iota_1;n)a_j}{k^G}I^*(\iota_1,\iota_2,j,k;n) +O(N^{-\delta})
\end{equation}
for some $\delta>0$, where we define the scaled modified Bessel function of the first kind 
\begin{equation}\label{eq:I_def}
    I^*(\iota_1,\iota_2,j,k;n) := \left[2\pi n-\frac{\pi\chi(S)}{12}\right]^{(G-1)/2}
    \left[-\frac{2\pi H}{k^2}-\frac{2\pi j}{k^2 \ell}\right]^{(1-G)/2}I_v(s),
\end{equation}
where $v := 1 - G$ and $s := 2\sqrt{\left[2\pi n-\frac{\pi\chi(S)}{12}\right]\left[-\frac{2\pi H}{k^2}-\frac{2\pi j}{k^2\ell}\right]}$.


Thus, by (\ref{eq:stratify_s}),  we see that taking $N\to\infty$ in \eqref{eq:exact_formula_S} yields an exact formula for $\xi_S(r,\ell,0,1;n)$, proving Theorem \ref{thm:exact_formulas}.
\end{proof}



\begin{remark}
Our use of $S(\iota_1, \iota_2, N; n)$ to provide an exact formula for $\xi_S(r, \ell, 0, 1; n)$ in \eqref{eq:stratify_s} hints at the relationship between the circle method and Poincar\'e series. In the method of Poincar\'e series, one averages functions over a group action to get modular forms with prescribed principle part at a given cusp and then uses Poisson summation to get exact formulas for the coefficients of their Fourier expansions in terms of an infinite sum of products of Kloosterman sums and Bessel functions. Our application of the circle method provides a concrete way of seeing how terms in this sum arise from the contribution of the integral near each cusp representative for the cusp under consideration.
\end{remark}

\begin{remark}
The proof of Theorem \ref{thm:exact_formulas} differs most from the proof of the corresponding statement in {\text {\rm \cite{reu_paper}}} in the establishment of the Kloosterman sum bound
\begin{equation} \label{eq:Ksum_bound} B_k(j,\ell,\iota_1,n) = O(n^{1/3}k^{2/3+\varepsilon}), \end{equation}
where  $h'$ is restricted to an interval $0 \leq \sigma_1 \leq h' < \sigma_2 \leq k$ and $\chi(S) = \sigma(S)$. The proof is a modification of the method originally proposed by Lehner in {\text {\rm \cite{Lehner}}}.
 \end{remark}

\subsection{Sketch of the proof of Theorem~\ref{thm:equidistribution}: Equidistribution}
Our main tool for proving Theorem \ref{thm:equidistribution} is an asymptotic formula for $\xi_S(r_1,\ell_1,r_2,\ell_2;n)$. We obtain this using Theorem \ref{thm:exact_formulas} by the same process as in the proof of Corollary 7.2 in \cite{reu_paper}.
This formula shows
that 
$\min\{{ \Lambda(}kj_1/\ell_1,kj_2/\ell_2)/k^2\}$ determines the dominant exponential term in the asymptotic description of  $ \xi(j_1,\ell_1,j_2,\ell_2;n)$, so that 
$ \xi(j_1,\ell_1,j_2,\ell_2;n)$ dominates $\xi(j_1',\ell_1,j_2',\ell_2;n)$ asymptotically if $\min\{{ \Lambda(}kj_1/\ell_1,kj_2/\ell_2)/k^2\} < \min\{{ \Lambda(}kj_1'/\ell_1,kj_2'/\ell_2)/k^2\}$.
If equality holds,
a sequence dominates asymptotically if its corresponding modular form is of strictly larger weight.
If these modular forms are of the same weight, then
 $$\xi(j_1,\ell_1,j_2,\ell_2;n) \sim \alpha\xi(j_1',\ell_1,j_2',\ell_2;n) $$ for some $\alpha >0$. 
Thus one can see by \ref{eq:roots} that if ${ \Lambda(}0,0) < { \Lambda(}j_1/\ell_1,j_2/\ell_2)$ for all $(j_1,j_2) \not \equiv (0,0)$, then we have $(\ell_1,\ell_2)$-equidistribution for $\mathcal{R} = \Z/ \ell_1 \Z \times \Z/\ell_2\Z$. This, along with manipulations of (\ref{eq:roots}), allows one to prove that $S$ has $(\ell_1,\ell_2)$-equidistribution in all five cases in Theorem \ref{thm:equidistribution}.

If $ \xi(0,\ell_1,0,\ell_2) = o(\xi(j_1,\ell_1,j_2,\ell_2))$ for some $(j_1,j_2) \not \equiv (0,0)$, then $S$ does not have $(\ell_1,\ell_2)$-equidistribution. To prove this claim, choose $(j_1,j_2)$ so that $\xi(j_1,\ell_1,j_2,\ell_2;n) \neq o( \xi(j_1',\ell_1,j_2',\ell_2;n))$ for all $(j_1',j_2')$. 
If equidistribution held, then we would have 
$$
C_S(0,1,0,1) \sim \alpha^* C_S(0,\ell_1,0,\ell_2)
$$
for some $\alpha^* >0$, which is false by our assumption that $ \xi(0,\ell_1,0,\ell_2) = o(\xi(j_1,\ell_1,j_2,\ell_2))$. This fact allows us to prove that $S$ does not have $(\ell_1,\ell_2)$-equidistribution in the cases not included in Theorem \ref{thm:equidistribution}. The cases with $\chi(S) <0$ require Theorem 15.1 in \cite{OurBook}.

The most difficult case is where $\min\{{ \Lambda(}j_1/\ell_1,j_2/\ell_2)\} = { \Lambda(}j_1/\ell_1,j_2/\ell_2) = { \Lambda(}0,0)$ for some $(j_1,j_2) \not \equiv (0,0)$, where one must prove that the weight of $Z_S(\zeta_{\ell_1}^{j_1},\ell_2^{j_2};\tau)$ is greater than that of $Z_S(1,1;\tau)$. To accomplish this task, one must first make use of the equality $h^{0,0} = 1$, to show that in this case we must have $\gcd(\ell_1,\ell_2) = 1$. One must prove that in this case $h^{1,0} >0$ and $\chi(S) + \sigma(S) >0$, and thus conclude that ${ \Lambda(}j_1/\ell_1,j_2/\ell_2)$ never obtains its minimum for $j_1 \equiv 0$, $j_2 \not \equiv 0$. It follows that $\ell_2 \neq 1$. Also, one can now produce a uniform description of the weight of all $Z_S(\zeta_{\ell_1}^{j_1},\zeta_{\ell_2}^{j_2};\tau)$ such that ${ \Lambda(}j_1/\ell_1,j_2/\ell_2) = { \Lambda(}0,0)$. If this weight is less than or equal to that of ${ \Lambda(}0,0)$, then for all $(x,y) \in [1/3,2/3] \times [2/5,3/5]$, we have ${ \Lambda(}x,y)< { \Lambda(}0,0)$. It follows that there is some $(j_1,j_2) \not \equiv (0,0)$ such that ${ \Lambda(}j_1/\ell_1,j_2/\ell_2) < { \Lambda(}0,0)$, contradicting our initial assumption. Therefore the weight of $Z_S(\zeta_{\ell_1}^{j_1},\zeta_{\ell_2}^{j_2};\tau)$ is greater than  that of ${ \Lambda(}0,0)$.

The final statement of Theorem \ref{thm:equidistribution} follows from the fact that ${ \Lambda(}0,0) - { \Lambda(}1/2,1/2) = h^{1,0}/2$.



\subsection{The case of $S=\C\mathbb{P}^2$}


We illustrate Theorem \ref{thm:exact_formulas} and Theorem \ref{thm:equidistribution} with numerics where $S = \C\mathbb{P}^2$. For the purposes of illustrating Theorem \ref{thm:exact_formulas}, we consider
\begin{equation*}
    Z_S(\zeta_3, -1;\tau) = 1 + 2 q + 4 q^2 + 7 q^3 + 12 q^4 + 20 q^5 + \cdots.
\end{equation*}
While Theorem~\ref{thm:exact_formulas} furnishes an infinite sum in $k$, in Tables \ref{table_2} and \ref{table_75} we approximate our exact formula by summing $k$ up to $N$, where $N$ is $2$ and $75$, respectively.

\renewcommand{\arraystretch}{1}
\begin{table}[h] 
\begin{center}
\scalebox{0.9}{\begin{tabular}{|c|c|c|c|c|c|}\hline 
     $n$ & $1$ &$2$ & $3$  & $4$ & $5$   \\\hline 
    $\xi_{2,S}(1,3,1,2)$ & $1.9374...$ & $3.8920...$ & $7.0204...$ & $12.1616...$ & $20.0159...$ \\\hline
\end{tabular}}
\end{center}
\caption{Approximate values in Theorem \ref{thm:exact_formulas}, $N=2$} \label{table_2}
\end{table}


\begin{table}[h] 
\small
\begin{center}
\scalebox{0.9}{\begin{tabular}{|c|c|c|c|c|c|}\hline 
     $n$ & $1$ &$2$ & $3$  & $4$ & $5$ \\\hline 
    $\xi_{75,S}(1,3,1,2)$ & $1.9989...$ & $4.0005...$ & $6.9995...$ & $12.0010...$ & $19.9995...$ \\\hline
\end{tabular}}
\end{center}
\caption{Approximate values in Theorem \ref{thm:exact_formulas}, $N=75
$} \label{table_75}
\end{table}

Tables \ref{table_equidistribution} and \ref{table_zeros} show the asymptotic equidistribution of the Hodge numbers of $\C\mathbb{P}^2$, which falls into case (4) in Theorem \ref{thm:equidistribution}. 
For this purpose, we define the following proportions
$$\Theta^{r_1,r_2}_{\ell_1, \ell_2, S}(n): = \frac{\gamma_S(r_1,\ell_1, r_2, \ell_2;n)}{\sum_{\substack{ j_1 \mod{\ell_1} \\ j_2 \mod{\ell_2}}} \gamma_S(j_1,\ell_1, j_2, \ell_2;n)}.$$ 
In Table \ref{table_equidistribution}, where $\gcd(\ell_1 = 3, \ell_2 = 2) = 1$, we see asymptotic equidistribution, while in Table \ref{table_zeros}, where $\gcd(\ell_1 = 2, \ell_2 = 4) = 2$, we get asymptotic equidistribution when $r_1 \equiv r_2 \ \mod{2}$ and $0$ otherwise. These numerics also suggest many underlying equalities that exist amongst the $\gamma_S(r_1, \ell_1, r_2, \ell_2)$ for different values of $(r_1, r_2)$ as a result of the symmetries of $Z_S(\zeta_{\ell_1}^{r_1}, \zeta_{\ell_2}^{r_2}; \tau)$.
\renewcommand{\arraystretch}{1.5}
\begin{table}[h] 
\small
\begin{center}
\scalebox{0.9}{\begin{tabular}{|c|c|c|c|c|c|}\hline 
     $n$&$5$ & $10$ &$15$ & $20$  &$25$  \\\hline 
    $\Theta^{0,0}_{3,2,S}(n)$  & $ 0.2222... $  &  $0.1886...     $  &  $0.1752...     $  & $0.1708...     $  &$0.1687...     $  \\\hline
    $\Theta^{0,1}_{3,2,S}(n)$  & $ 0.1111... $  &  $0.1446...     $  &  $0.1582...     $  & $0.1624...     $  &$0.1646...     $  \\\hline
    $\Theta^{1,0}_{3,2,S}(n)$  & $ 0.1296... $  &  $0.1571...     $  &  $0.1619...     $  & $0.1646...     $  &$0.1655...     $  \\\hline
    $\Theta^{1,1}_{3,2,S}(n)$  & $ 0.2037... $  &  $0.1761...     $  &  $0.1712...     $  & $0.1686...     $  &$0.1677...     $  \\\hline
    $\Theta^{2,0}_{3,2,S}(n)$  & $ 0.1296... $  &  $0.1571...     $  &  $0.1619...     $  & $0.1646...     $  &$0.1655...     $  \\\hline
    $\Theta^{2,1}_{3,2,S}(n)$  & $ 0.2037... $  &  $0.1761...     $  &  $0.1712...     $  & $0.1686...     $  &$0.1677...     $  \\\hline
\end{tabular}}
\end{center}
\caption{Comparative asymptotic properties of $\gamma_{S}(r_1, 3, r_2, 2; n)$ } \label{table_equidistribution}
\end{table}

\begin{table}[H] 
\small
\begin{center}
\scalebox{0.9}{\begin{tabular}{|c|c|c|c|c|c|}\hline 
     $n$&$5$ & $10$ &$15$ & $20$  &$25$  \\\hline 
    $\Theta^{0,0}_{2,4,S}(n)$  & $0.2592...$ & $ 0.2545... $  &  $0.2503...     $  &  $0.2503...     $  & $0.2500...     $  \\\hline
    $\Theta^{0,2}_{2,4,S}(n)$  & $ 0.2222... $  &  $0.2484...     $  &  $0.2488...     $  & $0.2498...     $  &$0.2498...     $  \\\hline
    $\Theta^{1,1}_{2,4,S}(n)$  & $ 0.2592... $  &  $0.2484...     $  &  $0.2503...     $  & $0.2498...     $  &$0.2500...     $  \\\hline
    $\Theta^{1,3}_{2,4,S}(n)$  & $ 0.2592... $  &  $0.2484...     $  &  $0.2503...     $  & $0.2498...     $  &$0.2500...     $  \\\hline
    $\Theta^{r_1,r_2}_{2,4,S}(n)$  & $ 0 $  &  $ 0     $  &  $0     $  & $0     $  &$0     $  \\\hline
\end{tabular}}
\end{center}
\caption{Comparative asymptotic properties of $\gamma_{S}(r_1, 2, r_2, 4; n)$ } \label{table_zeros}
\end{table}


\end{document}